\documentclass{amsart}

\usepackage[a4paper,hmargin={2.5cm,2.5cm},vmargin={3.5cm,3.5cm}]{geometry}

\usepackage[leqno]{amsmath}
\usepackage{amssymb,amsthm}
\usepackage[mathscr]{euscript}
\usepackage{bbm}
\usepackage{dsfont}
\usepackage[all,arc]{xy}
\usepackage{stmaryrd}
\DeclareMathAlphabet{\mathpzc}{OT1}{pzc}{m}{it}

\theoremstyle{plain}
\newtheorem*{theorem}{Theorem}
\newtheorem*{lemma}{Lemma}
\newtheorem*{proposition}{Proposition}
\newtheorem*{corollary}{Corollary}

\theoremstyle{definition}

\newtheorem*{examples}{Examples}

\theoremstyle{remark}

\newenvironment{eqcond}{\begin{enumerate}}{\end{enumerate}}

\newcommand{\Rw}{\Rightarrow}

\newcommand{\hrw}{\hookrightarrow}

\newcommand{\frp}{\mathfrak{p}}
\newcommand{\frq}{\mathfrak{q}}

\newcommand{\frw}{\mathfrak{w}}
\newcommand{\frv}{\mathfrak{v}}
\newcommand{\frx}{\mathfrak{x}}
\newcommand{\fry}{\mathfrak{y}}
\newcommand{\frz}{\mathfrak{z}}

\newcommand{\calI}{\mathcal{I}}
\newcommand{\calL}{\mathcal{L}}

\newcommand{\calP}{\mathcal{P}}
\newcommand{\calR}{\mathcal{R}}

\newcommand{\frX}{\mathfrak{X}}

\DeclareMathOperator{\ev}{ev}
\DeclareMathOperator{\Id}{Id}

\DeclareMathOperator{\ForgetSET}{O}
\DeclareMathOperator{\ForgetToV}{S}
\DeclareMathOperator{\ForgetToVAd}{A}
\DeclareMathOperator{\MFunctor}{M}

\DeclareMathOperator{\colim}{colim}
\DeclareMathOperator{\yoneda}{\mathpzc{y}}
\DeclareMathOperator{\Sup}{Sup}

\newcommand{\Psh}{\mathcal{P}}

\newcommand{\mate}[1]{\,^\ulcorner\! #1^\urcorner}

\newcommand{\fspstr}[2]{\llbracket #1,#2\rrbracket}

\newcommand{\catfont}[1]{\mathsf{#1}}

\newcommand{\V}{\catfont{V}}

\newcommand{\two}{\catfont{2}}

\newcommand{\quantale}{(\V,\otimes,k)}

\newcommand{\Pplus}{\catfont{P}_{\!\!{_+}}}

\newcommand{\SET}{\catfont{Set}}

\newcommand{\TOP}{\catfont{Top}}

\newcommand{\ORD}{\catfont{Ord}}

\newcommand{\Mat}[1]{#1\text{-}\catfont{Rel}}
\newcommand{\Mod}[1]{#1\text{-}\catfont{Mod}}

\newcommand{\Cat}[1]{#1\text{-}\catfont{Cat}}

\newcommand{\CatSep}[1]{#1\text{-}\catfont{Cat}_\mathrm{sep}}
\newcommand{\CatCompl}[1]{#1\text{-}\catfont{Cat}_\mathrm{cpl}}

\newcommand{\CocomplPhi}[1]{#1\text{-}\catfont{Cocont}^\Phi}

\newcommand{\CocomplSepPhi}[1]{#1\text{-}\catfont{Cocont}_\mathrm{sep}^\Phi}
\newcommand{\CocomplSepVar}[2]{#2\text{-}\catfont{Cocont}_\mathrm{sep}^{#1}}

\renewcommand{\to}{\longrightarrow}
\renewcommand{\mapsto}{\longmapsto}
\newcommand{\relto}{{\longrightarrow\hspace*{-2.8ex}{\mapstochar}\hspace*{2.6ex}}}

\newcommand{\kto}{\relbar\joinrel\rightharpoonup}
\newcommand{\krelto}{\,{\kto\hspace*{-2.5ex}{\mapstochar}\hspace*{2.6ex}}}
\newcommand{\kmodto}{\,{\kto\hspace*{-2.8ex}{\circ}\hspace*{1.3ex}}}

\newcommand{\kleisli}{\circ}

\newcommand{\multimapinv}{\mathrel{\circ\!-}}
\newcommand{\homkleislileft}{\multimap}
\newcommand{\homkleisliright}{\multimapinv}

\newcommand{\multimapdot}{\mathrel{-\!\bullet}}
\newcommand{\multimapdotinv}{\mathrel{\bullet\!-}}
\newcommand{\homcompleft}{\multimapdot}
\newcommand{\homcompright}{\multimapdotinv}

\newcommand{\Txi}{T_{\!_\xi}}

\newcommand{\monadfont}[1]{\mathbbm{#1}}
\newcommand{\mT}{\monadfont{T}}
\newcommand{\mI}{\monadfont{1}}
\newcommand{\mU}{\monadfont{U}}
\newcommand{\mL}{\monadfont{L}}
\newcommand{\mPhi}{\monadfont{I}^\Phi}

\newcommand{\monad}{(T,e,m)}

\newcommand{\imonad}{(\Id,1,1)}
\newcommand{\umonad}{(U,e,m)}
\newcommand{\wmonad}{(L,e,m)}
\newcommand{\Phimonad}{(\Phi,\yoneda^{\Phi},(\yoneda^{\Phi})^{-1})}

\newcommand{\theoryfont}[1]{\mathscr{#1}}

\newcommand{\Tth}{\theoryfont{T}}
\newcommand{\Ith}{\theoryfont{I}}
\newcommand{\Uth}{\theoryfont{U}}

\newcommand{\toptheory}{(\mT,\V,\xi)}
\newcommand{\itheory}{(\mI,\V,1_{\V})}

\newcommand{\sptheory}{(\mU,\two,\xi_{\two})}

\newcommand{\BC}{(BC)}

\newcommand{\op}{\mathrm{op}}
\newcommand{\co}{\mathrm{co}}

\usepackage[a4paper,hypertexnames=false]{hyperref} 

\begin{document}

\title[Relative injectivity as cocompleteness for a class of distributors]{Relative injectivity as cocompleteness\vspace*{1ex}\\ for a class of distributors}

\author{Maria Manuel Clementino}
\address{Departamento de Matem\'{a}tica\\ Universidade de Coimbra\\ 3001-454 Coimbra\\ Portugal}
\email{mmc@mat.uc.pt}

\author{Dirk Hofmann}
\address{Departamento de Matem\'{a}tica\\ Universidade de Aveiro\\3810-193 Aveiro\\ Portugal}
\email{dirk@ua.pt}

\thanks{The authors acknowledge partial financial assistance by Centro de Matem\'{a}tica da Universidade de Coimbra/FCT and Unidade de Investiga\c{c}\~{a}o e Desenvolvimento Matem\'{a}tica e Aplica\c{c}\~{o}es da Universidade de Aveiro/FCT}
\dedicatory{Dedicated to Walter Tholen on the occasion of his sixtieth birthday}

\maketitle

\begin{abstract}
Notions and techniques of enriched category theory can be used to study topological structures, like metric spaces, topological spaces and approach spaces, in the context of topological theories. Recently in [D. Hofmann, Injective spaces via adjunction, arXiv:math.CT/0804.0326] the construction of a Yoneda embedding allowed to identify injectivity of spaces as cocompleteness and to show monadicity of the category of injective spaces and left adjoints  over $\SET$. In this paper we generalise these results, studying cocompleteness with respect to a given class of distributors. We show in particular that the description of several semantic domains presented in [M. Escard\'o and B. Flagg, Semantic domains, injective spaces and monads,
Electronic Notes in Theoretical Computer Science 20 (1999)] can be translated into the $\V$-enriched setting.
\end{abstract}

\section*{Introduction}
This work continues the research line of previous papers, aiming to use categorical tools in the study of topological structures. Indeed, the perspective proposed in \cite{CH_TopFeat, CT_MultiCat} of looking at topological structures as (Eilenberg-Moore) lax algebras and, simultaneously, as a monad enrichment of $\V$-enriched categories, has shown to be very effective in the study of special morphisms -- like effective descent and exponentiable ones -- at a first step \cite{CH_EffDesc, CH_ExpVCat}, and recently in the study of (Lawvere/Cauchy-)completeness and injectivity \cite{CH_Compl, HT_LCls, Hof_Cocompl}. The results we present here complement this study of injectivity. More precisely, in the spirit of Kelly-Schmitt \cite{KS_Colim} we generalise the results of \cite{Hof_Cocompl}, showing that injectivity and cocompleteness -- when considered relative to a class of distributors -- still coincide. Suitable choices of this class of distributors allow us to recover, in the $\V$-enriched setting, results on injectivity of Escard\'{o}-Flagg \cite{EF_SemDom}.

The starting point of our study of injectivity is the notion of distributor (or bimodule, or profunctor), which allowed the study of weighted colimit, presheaf category, and the Yoneda embedding. It was then a natural step to `relativize' these ingredients and to consider {\em cocompleteness with respect to a class of distributors} $\Phi$. Namely, we introduce the notion of $\Phi$-cocomplete category, we construct the $\Phi$-presheaf category, and we prove that $\Phi$-cocompleteness is equivalent to the existence of a left adjoint of the Yoneda embedding into the $\Phi$-presheaf category. Furthermore, the class $\Phi$ determines a class of embeddings so that the injective $\Tth$-categories with respect to this class are precisely the $\Phi$-cocomplete categories. This result links our work with \cite{EF_SemDom}, where the authors study systematically semantic domains and injectivity characterisations with the help of Kock-Z\"{o}berlein monads.

\section{The Setting}
Throughout this paper we consider a (strict) \emph{topological theory} as introduced in \cite{Hof_TopTh}. Such a theory $\Tth=\toptheory$ consists of:
\begin{enumerate}
\item a commutative quantale $\V=\quantale$,
\item a $\SET$-monad $\mT=\monad$, where $T$ and $m$ satisfy \BC; that is, $T$ sends pullbacks to weak pullbacks and each naturality square of $m$ is a weak pullback, and
\item a $\mT$-algebra structure $\xi:T\V\to\V$ on $\V$ such that:
\begin{enumerate}
\item $\otimes:\V\times\V\to\V$ and $k:1\to\V$, $\ast\mapsto k$, are $\mT$-algebra homomorphisms making $(\V,\xi)$ a monoid in $\SET^\mT$; that is, the following diagrams
\begin{align*}
\xymatrix{T1\ar[d]_{!}\ar[r]^{Tk} & T\V\ar[d]^\xi\\ 1\ar[r]_k & \V} &&&
\xymatrix{T(\V\times\V)\ar[rr]^-{T(\otimes)}\ar[d]_{\langle\xi\cdot T\pi_1,\xi\cdot T\pi_2\rangle} && T\V\ar[d]^\xi\\ \V\times\V\ar[rr]_-{\otimes} && \V}
\end{align*}
are commutative;
\item For each set $X$, $\xi_X:\V^X\to \V^{TX}$, $(X\stackrel{\varphi}{\to}\V)\mapsto(TX\stackrel{T\varphi}{\to}T\V\stackrel{\xi}{\to}\V)$, defines a natural transformation $(\xi_{X})_{X}:P\to PT:\SET\to\ORD$.
\end{enumerate}
\end{enumerate}
Here $P:\SET\to\ORD$ is the $\V$-powerset functor defined as follows. We put $PX=\V^X$ with the pointwise order. Each map $f:X\to Y$ defines a monotone map $\V^f:\V^Y\to\V^X,\,\varphi\mapsto\varphi\cdot f$. Since $\V^f$ preserves all infima and all suprema, it has a left adjoint $Pf$. Explicitly, for $\varphi\in\V^X$ we have $Pf(\varphi)(y)=\bigvee\{\varphi(x)\mid x\in X,\,f(x)=y\}$.

\begin{examples}\label{ExTheories} Throughout this paper we will keep in mind the following topological theories:
\begin{enumerate}
\item The identity theory $\Ith=\itheory$, for each quantale $\V$, where $\mI=\imonad$ denotes the identity monad.
\item\label{ExTop} $\Uth_{\two}=\sptheory$, where $\mU=\umonad$ denotes the ultrafilter monad and $\xi_{\two}$ is essentially the identity map.
\item\label{ExApp} $\Uth_{\Pplus}=(\mU,\Pplus,\xi_{\Pplus})$ where $\Pplus=([0,\infty]^\op,+,0)$ and
\[
\xi_{\Pplus}:U\Pplus\to\Pplus,\;\;\frx\mapsto\inf\{v\in\Pplus\mid[0,v]\in\frx\}.
\]
\item\label{WordTh} The word theory $(\mL,\V,\xi_{_\otimes})$, for each quantale $\V$, where $\mL=\wmonad$ is the word monad and
\begin{align*}
\xi_{_\otimes}:L\V &\to\V.\\
(v_1,\ldots,v_n) &\mapsto v_1\otimes\ldots\otimes v_n\\
() &\mapsto k
\end{align*}
\end{enumerate}
\end{examples}

Every topological theory $\Tth=\toptheory$ encompasses several interesting ingredients.

\subsection*{\bf I} {\em The quantaloid $\Mat{\V}$} with {\em sets} as objects and {\em $\V$-relations} (also called {\em $\V$-matrices}, see \cite{BCRW_VarEnr}) $r:X\times Y\to\V$ as morphisms. We use the usual notation for relations, denoting the $\V$-relation $r:X\times Y\to\V$ by $r:X\relto Y$.
Since every map $f:X\to Y$ can be thought of as a $\V$-relation $f:X\times Y\to\V$ through its graph, there is an injective on objects and faithful functor $\SET\to\Mat{\V}$, unless $\V$ is degenerate (i.e. $k$ is the bottom element). Moreover, $\Mat{\V}$ has an involution $(-)^\circ:\Mat{\V}\to\Mat{\V}$, assigning to $r:X\relto Y$ the $\V$-relation $r^\circ:Y\relto X$, with $r^\circ(y,x):=r(x,y)$. For each $\V$-relation $r:X\relto Y$, the maps \[(-)\cdot r:\Mat{\V}(Y,Z)\to\Mat{\V}(X,Z)\mbox{ and }r\cdot (-):\Mat{\V}(Z,X)\to \Mat{\V}(Z,Y)\]
preserve suprema; hence they have right adjoints,
\[(-)\homcompright r:\Mat{\V}(X,Z)\to\Mat{\V}(Y,Z)\mbox{ and }r\homcompleft (-):\Mat{\V}(Z,Y)\to\Mat{\V}(Z,X).\]

\subsection*{\bf II} {\em The $\SET$-functor $T$ extends to a 2-functor} $\Txi:\Mat{\V}\to\Mat{\V}$ . To each $\V$-relation $r:X\times Y\to\V$, $\Txi$ assigns a $\V$-relation $\Txi r:TX\times TY \to\V$, which is the smallest (order-preserving) map $s:TX\times TY\to\V$ such that $\xi\cdot Tr\le s\cdot \langle T\pi_1,T\pi_2\rangle$.
\[
\xymatrix{T(X\times Y)\ar[rr]^{\langle T\pi_1,T\pi_2\rangle}\ar[dr]_{\xi_{X\times Y}(r)=\xi\cdot Tr} & &
TX\times TY\ar@{..>}[dl]^{\Txi r}\\ & \V\ar@{}[u]|-{\le}}
\]
Hence, for $\frx\in TX$ and $\fry\in TY$,
\[
\Txi r(\frx,\fry)=\bigvee\left\{\xi\cdot Tr(\frw)\;\Bigl\lvert\;\frw\in T(X\times Y), T\pi_1(\frw)=\frx,T\pi_2(\frw)=\fry\right\}.
\]
This 2-functor $\Txi$ preserves the involution, i.e. $\Txi(r^\circ)=\Txi(r)^\circ$ (and we write $\Txi r^\circ$) for each $\V$-relation $r:X\relto Y$, $m$ becomes a natural transformation $m:\Txi\Txi\to\Txi$ and $e$ an op-lax natural transformation $e:\Id\to\Txi$, i.e.\ $e_Y\circ r\leq \Txi r\circ e_X$ for all $r:X\relto Y$ in $\Mat{\V}$.

\subsection*{\bf III} A $\V$-relation of the form $\alpha:TX\relto Y$, called a \emph{$\Tth$-relation} and denoted by $\alpha:X\krelto Y$, will play an important role here. Given two $\Tth$-relations $\alpha:X\krelto Y$ and $\beta:Y\krelto Z$, their \emph{Kleisli convolution} $\beta\kleisli\alpha:X\krelto Z$ is defined as
\[\beta\kleisli\alpha=\beta\cdot \Txi\alpha\cdot m_X^\circ.\]
This operation is associative and has the $\Tth$-relation $e_X^\circ:X\krelto X$ as a lax identity: $a\kleisli e_X^\circ=a$ and $e_Y^\circ\kleisli a\ge a$ for any $a:X\krelto Y$.

\subsection*{\bf IV} $\Tth$-relations satisfying the usual unit and associativity categorical rules define $\Tth$-categories: a \emph{$\Tth$-category} is a pair $(X,a)$ consisting of a set $X$ and a $\Tth$-relation $a:X\krelto X$ on $X$ such that
\begin{align*}
e_X^\circ&\le a &&\text{and}& a\kleisli a&\le a.
\end{align*}
Expressed elementwise, these conditions become
\begin{align*}
k&\le a(e_X(x),x) &&\text{and}& \Txi a(\frX,\frx)\otimes a(\frx,x)\le a(m_X(\frX),x)
\end{align*}
for all $\frX\in TTX$, $\frx\in TX$ and $x\in X$. A function $f:X\to Y$ between $\Tth$-categories $(X,a)$ and $(Y,b)$ is a \emph{$\Tth$-functor} if $f\cdot a\le b\cdot Tf$, which in pointwise notation reads as
\[
 a(\frx,x)\le b(Tf(\frx),f(x))
\]
for all $\frx\in TX$, $x\in X$. The category of $\Tth$-categories and $\Tth$-functors is denoted by $\Cat{\Tth}$.

\subsection*{\bf V} In particular, the quantale $\V$ is a $\Tth$-category $\V=(\V,\hom_\xi)$, where \[\hom_\xi:T\V\times\V\to\V,\;(\frv,v)\mapsto\hom(\xi(\frv),v).\]

\subsection*{\bf VI} {\em The forgetful functor} $\ForgetSET:\Cat{\Tth}\to\SET,\;(X,a)\mapsto X$, \emph{is topological}, hence it has a left and a right adjoint. In particular, the free $\Tth$-category on a one-element set is given by $G=(1,e_1^\circ)$.

\subsection*{\bf VII} A $\V$-relation $\varphi:X\krelto Y$ between $\Tth$-categories $X=(X,a)$ and $Y=(Y,b)$ is a \emph{$\Tth$-distributor}, denoted as $\varphi:X\kmodto Y$, if $\varphi\kleisli a\le\varphi$ and $b\kleisli \varphi\le \varphi$. Note that we always have $\varphi\kleisli a\ge\varphi$ and $b\kleisli \varphi\ge \varphi$, so that the $\Tth$-distributor conditions above are in fact equalities. $\Tth$-categories and $\Tth$-distributors form a 2-category, denoted by $\Mod{\Tth}$, with Kleisli convolution as composition and with the 2-categorical structure inherited from $\Mat{\V}$.

\subsection*{\bf VIII} {\em Each $\Tth$-functor $f:(X,a)\to(Y,b)$ induces an adjunction $f_*\dashv f^*$ in $\Mod{\Tth}$}, with $f_*:X\krelto Y$ and $f^*:Y\krelto X$ defined as $f_*=b\cdot Tf$ and $f^*=f^\circ\cdot b$ respectively. In fact, these assignments are functorial, i.e. they define two functors:
\[\begin{array}{rclcrcl}
 (-)_*:\Cat{\Tth}^\co&\to&\Mod{\Tth} &\text{and}& (-)^*:\Cat{\Tth}^\op&\to&\Mod{\Tth},\\
 X&\mapsto&X_*=X&&X&\mapsto&X^*=X\\
 f&\mapsto&f_*=b\cdot Tf&&f&\mapsto&f^*=f^\circ\cdot b
\end{array}\]
A $\Tth$-functor $f:X\to Y$ is called \emph{fully faithful} if $f^*\kleisli f_*=1_X^*$, while it is called \emph{dense} if $f_*\kleisli f^*=1_Y^*$. Note that $f$ is fully faithful if and only if, for all $\frx\in TX$ and $x\in X$, $a(\frx,x)=b(Tf(\frx),f(x))$.

\subsection*{\bf IX} For a $\Tth$-distributor $\alpha:X\kmodto Y$, {\em the composition function $-\kleisli\alpha$ has a right adjoint} $(-)\homkleisliright\alpha$ where, for a given $\Tth$-distributor $\gamma:X\kmodto Z$, the extension $\gamma\homkleisliright\alpha:Y\kmodto Z$ is constructed in $\Mat{\V}$ as the extension $\gamma\homkleisliright \alpha=\gamma\homcompright (\Txi\alpha\cdot m_X^\circ)$.
\[
\xymatrix{TX\ar[r]|-{\object@{|}}^\gamma\ar[d]|-{\object@{|}}_{m_X^\circ} & Z.\\
TTX\ar[d]|-{\object@{|}}_{\Txi\alpha}\\ TY\ar@{.>}[ruu]|-{\object@{|}}}
\]
The following rules are easily checked.
\begin{lemma}\label{CalRulesMod}
The following assertions hold.
\begin{enumerate}
\item If $\alpha$ is a right adjoint, then $\alpha\kleisli(\varphi\homkleisliright\psi)= (\alpha\kleisli\varphi)\homkleisliright\psi$.
\item If $\gamma\dashv\delta$, then $(\alpha\homkleisliright\beta)\kleisli\gamma=\alpha\homkleisliright(\delta\kleisli\beta)$.
\item If $\gamma\dashv\delta$, then $(\alpha\kleisli\gamma)\homkleisliright\beta=\alpha\homkleisliright(\beta\kleisli\delta)$.
\end{enumerate}
\end{lemma}

\subsection*{\bf X}
It is also important the interplay of several functors relating the structures, i.e. {\em Eilenberg-Moore algebras}, {\em $\Tth$-categories} and {\em $\V$-categories}.
The inclusion functor $\SET^\mT\hrw\Cat{\Tth}$, given by regarding the structure map $\alpha:TX\to X$ of an Eilenberg-Moore algebra $(X,\alpha)$ as a $\Tth$-relation $\alpha:X\krelto X$, has a left adjoint, constructed \emph{\`{a} la \v Cech-Stone compactification} in \cite{CH_TopFeat}.
 \begin{align*}
\xymatrix{\SET^\mT\ar@<-1ex>@{^(->}[rr]^\bot && \Cat{\Tth}\ar@/_1pc/[ll]}
\end{align*}
We denote by $|X|$ the free Eilenberg-Moore algebra $(TX,m_X)$ considered as a $\Tth$-category.

Making use of the identity $e:\Id\to T$ of the monad, to each $\Tth$-category $X=(X,a)$ we assign a $\V$-category structure on $X$, $a\cdot e_X:X\relto X$. This correspondence defines a functor $\ForgetToV:\Cat{\Tth}\to\Cat{\V}$, which has also a left adjoint $\ForgetToVAd:\Cat{\V}\to\Cat{\Tth}$, with $\ForgetToVAd (X,a):=(X,e_X^\circ\cdot\Txi r)$.
 \begin{align*}
\xymatrix{\Cat{\Tth}\ar@<-1ex>[rr]^\bot_\ForgetToV && \Cat{\V}.\ar@/_1pc/[ll]_\ForgetToVAd}
\end{align*}

Furthermore, making now use of the multiplication $m:T^2\to T$ of the monad, one can define a functor
\[
\MFunctor:\Cat{\Tth}\to\Cat{\V}
\]
which sends a $\Tth$-category $(X,a)$ to the $\V$-category $(TX,\Txi a\cdot m_X^\circ)$.

We can now define the process of \emph{dualizing a $\Tth$-category} as the composition of the following functors
\[\xymatrix{\Cat{\Tth}\ar[d]_M\ar@{..>}[r]^{(\;)^\op}&\Cat{\Tth}\\
\Cat{\V}\ar[r]_{(\;)^\op}&\Cat{\V}\ar[u]_A}\]
that is, the \emph{dual of a $\Tth$-category $(X,a)$} is defined as
\[X^\op=\ForgetToVAd(\MFunctor(X)^\op),\]
which is a structure on $TX$. If $\mT$ is the identity monad, then $X^\op$ is indeed the dual $\V$-category of $X$.

\subsection*{\bf XI}
The tensor product on $\V$ can be transported to $\Cat{\Tth}$ by putting \[(X,a)\otimes(Y,b)=(X\times Y,c),\] with
\[
c(\frw,(x,y))=a(T\pi_1(\frw),x)\otimes b(T\pi_2(\frw),y),
\]
where $\frw\in T(X\times Y)$, $x\in X$, $y\in Y$. The $\Tth$-category $E=(1,k)$ is a $\otimes$-neutral object, where $1$ is a singleton set and $k:T1\times 1\to\V$ the constant relation with value $k\in\V$. For each set $X$, the functor $|X|\otimes(-):\Cat{\Tth}\to\Cat{\Tth}$ has a right adjoint $(-)^{|X|}:\Cat{\Tth}\to\Cat{\Tth}$. Explicitly, the structure $\fspstr{-}{-}$ on $\V^{|X|}$ is given by the formula
\[
\fspstr{\frp}{\psi}=
\bigwedge_{\substack{\frq\in T(|X|\times\V^{|X|})\\ \frq\mapsto\frp}}\hom(\xi\cdot T\ev(\frq),\psi(m_X\cdot T\pi_1(\frq))),
\]
for each $\frp\in T\V^{|X|}$ and $\psi\in\V^{|X|}$.
\begin{theorem} {\em \cite{CH_Compl}}\label{CharTMod}
For $\Tth$-categories $(X,a)$ and $(Y,b)$, and a $\Tth$-relation $\psi:X\krelto Y$, the following assertions are equivalent.
\begin{eqcond}
\item $\psi:(X,a)\kmodto(Y,b)$ is a $\Tth$-distributor.
\item Both $\psi:|X|\otimes Y\to\V$ and $\psi:X^\op\otimes Y\to\V$ are $\Tth$-functors.
\end{eqcond}
\end{theorem}

\subsection*{\bf XII}
Hence, {\em each $\Tth$-distributor $\varphi:X\kmodto Y$ provides a $\Tth$-functor}
\[
\mate{\varphi}:Y\to\V^{|X|}
\]
which factors through the embedding $\Psh X\hrw\V^{|X|}$, where $\Psh X=\{\psi\in\V^{|X|}\mid \psi:X\kmodto G\}$ is the $\Tth$-category of \emph{contravariant presheafs on $X$}:
\[
\xymatrix{Y\ar[r]^{\mate{\varphi}}\ar[dr]_ {\mate{\varphi}}& \V^{|X|}\\ & \Psh X\ar@{_(->}[u]}
\]
In particular, for each $\Tth$-category $X=(X,a)$, the $\V$-relation $a:TX\times X\to\V$ is a $\Tth$-distributor $a:X\kmodto X$, and therefore we have the \emph{Yoneda functor}
\[
 \yoneda_X=\mate{a}:X\to\Psh X.
\]
\begin{theorem}{\em \cite{Hof_Cocompl}}\label{Yoneda}
Let $\psi:X\kmodto Z$ and $\varphi:X\kmodto Y$ be $\Tth$-distributors. Then, for all $\frz\in TZ$ and $y\in Y$,
\[
\fspstr{T\mate{\psi}(\frz)}{\mate{\varphi}(y)}=(\varphi\homkleisliright\psi)(\frz,y).
\]
\end{theorem}
\begin{corollary}{\em \cite{Hof_Cocompl}}\label{YonedaLem}
For each $\varphi\in\hat{X}$ and each $\frx\in TX$, $\varphi(\frx)=\fspstr{T\yoneda_X(\frx)}{\varphi}$, that is, $(\yoneda_X)_*:X\kmodto\hat{X}$ is given by the evaluation map $\ev:TX\times\hat{X}\to\V$. As a consequence, $\yoneda_X:X\to\hat{X}$ is fully faithful.
\end{corollary}

\subsection*{\bf XIII}
Transporting the order-structure on hom-sets from $\Mod{\Tth}$ to $\Cat{\Tth}$ via the functor $(-)^*:\Cat{\Tth}^\op\to\Mod{\Tth}$, {\em $\Cat{\Tth}$ becomes a 2-category}. That is, for $\Tth$-functors $f,g:X\to Y$ we define $f\le g$ whenever $f^*\le g^*$, which in turn is equivalent to $g_*\le f_*$. We call $f,g:X\to Y$ \emph{equivalent}, and write $f\cong g$, if $f\le g$ and $g\le f$. Hence, $f\cong g$ if and only if $f^*=g^*$ if and only if $f_*=g_*$. A $\Tth$-category $X$ is called \emph{separated} (see \cite{HT_LCls} for details) whenever $f\cong g$ implies $f=g$, for all $\Tth$-functors $f,g:Y\to X$ with codomain $X$. One easily verifies that the {\em $\Tth$-category $\V=(\V,\hom_\xi)$ is separated}, and so is each $\Tth$-category of the form $\Psh X$ for a $\Tth$-category $X$. The full subcategory of $\Cat{\Tth}$ consisting of all separated $\Tth$-categories is denoted by $\CatSep{\Tth}$. The 2-categorical structure on $\Cat{\Tth}$ allows us to consider adjoint $\Tth$-functors: $\Tth$-functor $f:X\to Y$ is \emph{left adjoint} if there exists a $\Tth$-functor $g:Y\to X$ such that $1_X\le g\cdot f$ and $1_Y\ge f\cdot g$. Considering the corresponding $\Tth$-distributors, $f$ is left adjoint to $g$ if and only if $g_*\dashv f_*$, that is, if and only if $f_*=g^*$.

A more complete study of this subject can be found in \cite{Hof_TopTh, Hof_Cocompl}.

\section{The results}\label{results}

In the sequel we consider a class $\Phi$ of $\Tth$-distributors subject to the following axioms.
\begin{description}
\item[(Ax~1).] For each $\Tth$-functor $f$, $f^*\in\Phi$.
\item[(Ax~2).] For all $\varphi\in\Phi$ and all $\Tth$-functors $f:A\to X$ we have
\begin{align*}
f^*\kleisli\varphi&\in\Phi, &\varphi\kleisli f^*&\in\Phi, & f_*\in\Phi \Rw \varphi\circ f_*\in\Phi;
\end{align*}
whenever the compositions are defined.
\item[(Ax~3).] For all $\varphi:X\kmodto Y\in\Mod{\Tth}$,
\[
(\forall y\in Y\,.\,y^*\kleisli \varphi\in\Phi)\Rw \varphi\in\Phi
\]
where $y^*$ is induced by $y:1\to Y$, $*\mapsto y$.
\end{description}
Condition (Ax~2) requires that $\Phi$ is closed under certain compositions. In fact, in most examples $\Phi$ will be closed under arbitrary compositions. Furthermore, there is a largest and a smallest such class of $\Tth$-distributors, namely the class $\calP$ of all $\Tth$-distributors and the class $\calR=\{f^*\mid f:X\to Y\}$ of all representable $\Tth$-distributors.

We call a $\Tth$-functor $f:X\to Y$ \emph{$\Phi$-dense} if $f_*\in\Phi$. Certainly, if $f$ is a left adjoint $\Tth$-functor, with $f\dashv g$, then $f_*=g^*\in\Phi$, i.e. $f$ is $\Phi$-dense. A $\Tth$-category $X$ is called \emph{$\Phi$-injective} if, for all $\Tth$-functors $f:A\to X$ and fully faithful $\Phi$-dense $\Tth$-functors $i:A\to B$, there exists a $\Tth$-functor $g:B\to X$ such that $g\cdot i\cong f$. Furthermore, $X$ is called \emph{$\Phi$-cocomplete} if each weighted diagram
\[
\xymatrix{Y\ar[r]^{h}\ar@{-^{>}}|-{\object@{o}}[d]_\varphi & X\\ Z}
\]
with $\varphi\in\Phi$ has a colimit $g\cong\colim(\varphi,h):Z\to X$. A $\Tth$-functor $f:X\to Y$ is \emph{$\Phi$-cocontinuous} if $f$ preserves all existing $\Phi$-weighted colimits. Note that in both cases it is enough to consider diagrams where $h=1_X$. We denote by $\CocomplPhi{\Tth}$ the 2-category of all $\Phi$-cocomplete $\Tth$-categories and $\Phi$-cocontinuous $\Tth$-functors, and by $\CocomplSepPhi{\Tth}$ its full subcategory of all $\Phi$-cocomplete and separated $\Tth$-categories.

If $\Phi$ is the class $\calP$ of all $\Tth$-distributors, then $\CocomplPhi{\Tth}$ is the category of cocomplete $\Tth$-categories and left adjoint $\Tth$-functors (as shown in \cite[Prop. 2.12]{Hof_Cocompl}).

\begin{lemma}\label{CompCancDense}
Consider the (up to isomorphism) commutative triangle
\[
\xymatrix{X\ar[d]_f\ar[dr]^h_\cong \\ Y\ar[r]_g & Z}
\]
of $\Tth$-functors. Then the following assertions hold.
\begin{enumerate}
\item If $g$ and $f$ are $\Phi$-dense, then so is $h$.
\item If $h$ is $\Phi$-dense and $g$ is fully faithful, then $f$ is $\Phi$-dense.
\item If $h$ is $\Phi$-dense and $f$ is dense, then $g$ is $\Phi$-dense.
\end{enumerate}
\end{lemma}
\begin{proof}
The proof is straightforward:
(1) $h_*=g_*\circ f_*\in\Phi$ by (Ax~2), since $g_*, f_*\in\Phi$;
(2) $f_*=g^*\circ g_*\circ f_*=g^*\circ h_*\in\Phi$ by (Ax~2), since $h_*\in\Phi$; (3) $g_*=g_*\circ f_*\circ f^*=h_*\circ f^*\in\Phi$ by (Ax~2), since $h_*\in\Phi$.
\end{proof}
We put now
\[
 \Phi X=\{\psi\in \Psh X\mid \psi\in\Phi\}
\]
considered as a subcategory of $\calP X$. We have the restriction
\[
\yoneda_X^\Phi:X\to\Phi X
\]
of the Yoneda map, and each $\psi\in\Phi X$ is a $\Phi$-weighted colimit of representables (see \cite[Proposition 2.5]{Hof_Cocompl}).
\begin{lemma}
The following assertions hold.
\begin{enumerate}
\item $\yoneda_X^\Phi:X\to\Phi X$ is $\Phi$-dense.
\item For each $\Tth$-distributor $\varphi:X\kmodto Y$, $\varphi\in\Phi$ if and only if $\mate{\varphi}:Y\to\Psh X$ factors through the embedding $\Phi X\hrw\Psh X$.
\end{enumerate}
\end{lemma}
\begin{proof}
By the Yoneda Lemma (Corollary~\ref{YonedaLem}), for any $\psi\in\Phi X$ we have $\psi^*\kleisli(\yoneda_X^\Phi)_*=\psi\in\Phi$, therefore $(\yoneda_X^\Phi)_*\in\Phi$ by (Ax~3) and the assertion (1) follows. To see (2), just observe that $\mate{\varphi}(y)=y^*\kleisli\varphi$, and use again (Ax~3).
\end{proof}
Our next result extends Theorem 2.6 of \cite{Hof_Cocompl}. We omit its proof because it uses exactly the same arguments.
\begin{theorem}\label{CharCocomplInj}
The following assertions are equivalent, for a $\Tth$-category $X$.
\begin{eqcond}
\item $X$ is $\Phi$-injective.
\item $\yoneda_X^\Phi:X\to\Phi X$ has a left inverse $\Sup_X^\Phi:\Phi X\to X$.
\item $\yoneda_X^\Phi:X\to\Phi X$ has a left adjoint $\Sup_X^\Phi:\Phi X\to X$.
\item $X$ is $\Phi$-cocomplete.
\end{eqcond}
\end{theorem}

Recall from \cite{Hof_Cocompl} that, for a given $\Tth$-functor $f:X\to Y$, we have an adjoint pair of $\Tth$-functors $\Psh f\dashv f^{-1}$ where
\begin{align*}
\Psh f: \Psh X &\to \Psh Y &\text{and}&& f^{-1}: \Psh Y &\to \Psh X.\\
\psi &\mapsto \psi\kleisli f^* &&& \psi &\mapsto \psi\kleisli f_*
\end{align*}
By (Ax~1) and (Ax~2), the $\Tth$-functor $\Psh f: \Psh X \to \Psh Y$ restricts to a $\Tth$-functor $\Phi f:\Phi X\to\Phi Y$. On the other hand, $f^{-1}: \Psh Y \to \Psh X$ restricts to $f^{-1}: \Phi Y \to \Phi X$ provided that $f$ is $\Phi$-dense.
\begin{proposition}
The following conditions are equivalent for a $\Tth$-functor $f:X\to Y$.
\begin{eqcond}
\item $f$ is $\Phi$-dense.
\item $\Phi f$ is left adjoint.
\item $\Phi f$ is $\Phi$-dense.\newcounter{counter}\setcounter{counter}{\value{enumi}}
\end{eqcond}
\end{proposition}
\begin{proof} (i) $\Rightarrow$ (ii): If $f$ is $\Phi$-dense, then $\Phi f\dashv f^{-1}:\Phi Y\to\Phi X$ defined above.
(ii) $\Rightarrow$ (iii): If $\Phi f\dashv g$, then $(\Phi f)_*=g^*\in\Phi$, i.e. $\Phi f$ is $\Phi$-dense.
(iii) $\Rightarrow$ (i): Consider the diagram
\[\xymatrix{X\ar[r]^{\yoneda_X^\Phi}\ar[d]_f&\Phi X\ar[d]^{\Phi f}\\
Y\ar[r]_{\yoneda_Y^\Phi}&\Phi Y}\]
If $\Phi f$ is $\Phi$-dense, then $\yoneda^\Phi_Y\cdot f=\Phi f\cdot \yoneda_X^\Phi$ is $\Phi$-dense, and so by \ref{CompCancDense}(2) $f$ is $\Phi$-dense because $\yoneda_Y^\Phi$ is fully faithful.
\end{proof}
In particular, for each $\Tth$-category $X$, $\Phi \yoneda_X^\Phi:\Phi X\to\Phi\Phi X$ has a right adjoint, $(\yoneda_X^\Phi)^{-1}$.
We show next that $(\yoneda_X^\Phi)^{-1}$ has also a right adjoint, $\yoneda_{\Phi X}^\Phi:\Phi X\to\Phi\Phi X$, so that:
\[\Phi \yoneda_X^\Phi\dashv (\yoneda_X^\Phi)^{-1}=\Sup_{\Phi X}^\Phi\dashv \yoneda_{\Phi X}^\Phi.\]
\begin{proposition}\label{PhiXcocompl}
For each $\Tth$-category $X$, $\Phi X$ is $\Phi$-cocomplete where $\Sup_{\Phi X}^\Phi=(\yoneda_X^\Phi)^{-1}$.
\end{proposition}
\begin{proof}
Since $\yoneda_X^\Phi$ is $\Phi$-dense, we may define $\Sup_{\Phi X}^\Phi:=(\yoneda_X^\Phi)^{-1}$. We have to show that $\Sup_{\Phi X}^\Phi$ is a left inverse for $\yoneda_{\Phi X}^\Phi$; that is, $(\yoneda_X^\Phi)^{-1}\cdot \yoneda_{\Phi X}^\Phi=1_{\Phi X}$: for each $\psi\in\Phi X$, $((\yoneda_X^\Phi)^{-1}\cdot\yoneda_{\Phi X}^\Phi)(\psi)=\psi^*\circ (\yoneda_X^\Phi)_*=\psi$.
\end{proof}
In \cite{Hof_Cocompl} we constructed $\Psh f$ as the colimit $\Psh f\cong\colim((\yoneda_X)_*,\yoneda_Y\cdot f)$, and a straightforward calculation shows that also $\Phi f\cong\colim((\yoneda_X^\Phi)_*,\yoneda_Y^\Phi\cdot f)$, for each $\Tth$-functor $f:X\to Y$. To see this, we consider the commutative diagrams
\[
\xymatrix{X\ar[r]^{\yoneda_X^\Phi}\ar[d]_f\ar@/^1.8pc/[rr]^{\yoneda_X} &
\Phi X,\ar[r]^{i_X}\ar[d]^{\Phi f} & \Psh X\ar[d]^{\Psh f}\\
Y\ar[r]_{\yoneda_Y^\Phi}\ar@/_1.8pc/[rr]_{\yoneda_Y} & \Phi
Y\ar[r]_{i_Y} & \Psh Y}
\]
and obtain
\begin{align*}
(\Phi f)_* &= i_Y^*\kleisli {i_Y}_*\kleisli(\Phi f)_*\\
&= i_Y^*\kleisli(\Psh f)_*\kleisli {i_X}_*\\
&=i_Y^*\kleisli(({\yoneda_Y}_*\kleisli
f_*)\homkleisliright{\yoneda_X}_*)\kleisli{i_X}_* &&\text{since $\Psh
f\cong\colim((\yoneda_X)_*,\yoneda_Y\cdot f)$}\\
&= (i_Y^*\kleisli{\yoneda_Y}_*\kleisli
f_*)\homkleisliright({i_X}^*\kleisli{\yoneda_X}_*) &&\text{by Lemma
\ref{CalRulesMod}}\\
&= ({\yoneda^\Phi_Y}_*\kleisli f_*)\homkleisliright{\yoneda^\Phi_X}_*.
\end{align*}

\begin{proposition}\label{CharPhiCont}
Let $f:X\to Y$ a $\Tth$-functor where $X$ and $Y$ are $\Phi$-cocomplete.
\begin{enumerate}
\item The following assertions are equivalent.
\begin{eqcond}
\item $f$ is $\Phi$-cocontinuous.
\item We have $f\cdot \Sup^\Phi_X\cong \Sup^\Phi_Y\cdot\Phi f$.
\[
 \xymatrix{\Phi X\ar[r]^{\Phi f}\ar[d]_{\Sup^\Phi_X}\ar@{}[dr]|{\cong}
	& \Phi Y\ar[d]^{\Sup^\Phi_Y}\\ X\ar[r]_f & Y}
\]
\end{eqcond}
\item If $f$ is $\Phi$-cocontinuous, then $f$ is $\Phi$-dense if and only it is a left adjoint.
\end{enumerate}
\end{proposition}
\begin{proof} (1) (a) $\Rightarrow$ (b): Recall that
\[\xymatrix{X\ar[r]^{1_X}\ar@{-^{>}}|-{\object@{o}}[d]_{(\yoneda_X^\Phi)_*} & X\\
\Phi X\ar@{.>}[ru]_{(\Sup_X^\Phi)_*=1_X\homkleisliright (\yoneda_X^\Phi)_*}}\]
Hence
\[\begin{array}{rcl}
(f\cdot \Sup_X^\Phi)_*&=&f_*\homkleisliright (\yoneda_X^\Phi)_*\\
&=&((\Sup_Y^\Phi)_*\circ (\yoneda_Y^\Phi)_*\circ f_*)\homkleisliright (\yoneda_X^\Phi)_*\\
&=&(\Sup_Y^\Phi)_*\circ((\yoneda_Y^\Phi)_*\circ f_*\homkleisliright (\yoneda_X^\Phi)_*)\\
&=&(\Sup_Y^\Phi)_*\circ \Phi f_*.
\end{array}\]

(b)$\Rightarrow$ (a): Consider
\[\xymatrix{X\ar@{-^{>}}|-{\object@{o}}[r]^{1_X^*}\ar@{-^{>}}|-{\object@{o}}[d]_\varphi & X\ar[r]^f&Y\\
A\ar@{.>}[ru]_{(\Sup_X^\Phi\cdot\mate{\varphi})_*}&&}\]
Then
\[\begin{array}{rcl}
(f\cdot \Sup_X^\Phi\cdot \mate{\varphi}) &=&\Sup_Y^\Phi\cdot \Phi f\cdot \mate{\varphi}\\
&=&\Sup_Y^\Phi\cdot \mate{\varphi\cdot {f^*} }\\
&\cong& \colim(\varphi, f)
\end{array}\]

(2) If $f$ is $\Phi$-cocontinuous and $\Phi$-dense, from the commutative diagram of (1)(b) we have
$f\dashv \Sup_X^\Phi\cdot f^{-1}\cdot \yoneda_Y^\Phi$ since $f\cdot \Sup_X^\Phi=\Sup_Y^\Phi\cdot\Phi f\dashv f^{-1}\cdot \yoneda_Y^\Phi$ and $\Sup_X^\Phi\cdot\yoneda_X^\Phi=1_X$. The converse is trivially true.
\end{proof}
\begin{corollary}
$\Phi X$ is closed in $\Psh X$ under $\Phi$-weighted colimits.
\end{corollary}
\begin{proof}
We show that the inclusion functor $i:\Phi X\to\Psh X$ is $\Phi$-cocontinuous, which, by the proposition above, is equivalent to the commutativity of the diagram
\[
\xymatrix{\Phi \Phi X\ar[r]^{\Phi i}\ar[d]_{\Sup^\Phi_{\Phi X}} &
\Phi \Psh X\ar[d]^{\Sup^\Phi_{\Psh X}}\\ \Phi X\ar[r]_i & \Psh X.}
\]
In Proposition \ref{PhiXcocompl} we observed $\Sup^\Phi_{\Phi X}=(\yoneda_X^\Phi)^{-1}$, and from Theorem \ref{CharCocomplInj} and \cite[Theorem~2.8]{Hof_Cocompl} follows that $\Sup^\Phi_{\Psh X}$ is the restriction of $\yoneda_X^{-1}:\Psh\Psh X\to\Psh X$ to $\Phi\Psh X$. Let $\Psi\in\Phi\Phi X$. Then
\begin{align*}
i\cdot(\yoneda_X^\Phi)^{-1}(\Psi)&=\Psi\kleisli(\yoneda_X^\Phi)_*
\intertext{and}
\yoneda_X^{-1}\cdot\Phi i(\Psi) &= \yoneda_X^{-1}(\Psi\kleisli i^*)=\Psi\kleisli i^*\kleisli(\yoneda_X)_*=\Psi\kleisli(\yoneda_X^\Phi)_*,
\end{align*}
and the assertion follows.
\end{proof}
Theorem \ref{CharCocomplInj} says in particular that, for each $\Tth$-functor $f:A\to X$, $\Phi$-injective $\Tth$-category $X$ and fully faithful $\Phi$-dense $\Tth$-functor $i:A\to B$, we have a canonical extension $g:B\to X$ of $f$ along $i$, namely $g\cong\colim(i_*,f)$, giving us an alternative description of $\Phi f$.
\begin{theorem}\label{Kan}
Composition with $\yoneda_X^\Phi:X\to\Phi X$ defines an equivalence
\[
 \CocomplPhi{\Tth}(\Phi X,Y)\to\Cat{\Tth}(X,Y)
\]
of ordered sets, for each $\Phi$-cocomplete $\Tth$-category $Y$.
\end{theorem}

The series of results above tell us that $\CocomplSepPhi{\Tth}$ is actually a (non-full) reflective subcategory of $\Cat{\Tth}$, with left adjoint $\Phi:\Cat{\Tth}\to\CocomplSepPhi{\Tth}$. In fact, $\Phi$ is a 2-functor and one verifies as in \cite{Hof_Cocompl} that the induced monad $\mPhi=\Phimonad$ on $\Cat{\Tth}$ is of Kock-Z\"{o}berlein type. Theorem~\ref{CharCocomplInj} and Proposition~\ref{CharPhiCont} imply that $\CocomplSepPhi{\Tth}$ is equivalent to the category of Eilenberg-Moore algebras of $\mPhi$.

Finally, we wish to study monadicity of the canonical forgetful functor \[G:\CocomplSepPhi{\Tth}\to\SET.\]
Certainly,
\begin{description}
\item[(a)] $G$ has a left adjoint given by the composite
\[
 \SET\xrightarrow{\hspace{1em}\text{disc}\hspace{1em}}\Cat{\Tth}
\xrightarrow{\hspace{1em}\Phi\hspace{1em}}\CocomplSepPhi{\Tth},
\]
where $\text{disc}(X)=(X,e_X^\circ)$, and $\text{disc}(f)=f$. 
\end{description}
In order to prove monadicity of $G$ we will impose, in addition to (Ax~1)-(Ax~3),
\begin{description}
\item[(Ax~4).] For each surjective $\Tth$-functor $f$, $f_*\in\Phi$.
\end{description}
Hence, any bijective $f:X\to Y$ in $\CocomplSepPhi{\Tth}$ is $\Phi$-dense and therefore left adjoint. By \cite[Lemma 2.16]{Hof_Cocompl}, $f$ is invertible and we have seen that 
\begin{description}
\item[(b)] $G$ reflects isomorphisms. 
\end{description}
In order to conclude that $G$ is monadic, it is left to show that 
\begin{description}
\item[(c)] $\CocomplSepPhi{\Tth}$ has and $G$ preserves coequaliser of $G$-equivalence relations
 \end{description}
(see, for instance, \cite[Corollary 2.7]{MS_Monads}). To do so, let $\pi_1,\pi_2:R\rightrightarrows X$ in $\CocomplSepPhi{\Tth}$ be an equivalence relation in $\SET$, where $\pi_1$ and $\pi_2$ are the projection maps, and let $q:X\to Q$ be its coequaliser in $\Cat{\Tth}$. The proof in \cite[Section 2.6]{Hof_Cocompl} rests on the observation that
\begin{equation*}
\xymatrix{\Psh R\ar@<0.5ex>[r]^{\Psh\pi_1}\ar@<-0.5ex>[r]_{\Psh\pi_2} &\Psh X\ar[r]^{\Psh q} &\Psh Q}
\end{equation*}
is a split fork in $\CatSep{\Tth}$. Naturally, we wish to show that, in our setting,
\begin{equation*}
\xymatrix{\Phi R\ar@<0.5ex>[r]^{\Phi\pi_1}\ar@<-0.5ex>[r]_{\Phi\pi_2} &\Phi X\ar[r]^{\Phi q} &\Phi Q}
\end{equation*}
gives rise to a split fork in $\CatSep{\Tth}$ as well. Since $\pi_1$, $\pi_2$ and $q$ are surjective, the $\Tth$-functors $\pi_1$, $\pi_2$ and $q$ are $\Phi$-dense and therefore we have $\Tth$-functors $q^{-1}:\Phi Q\to\Phi X$ and $\pi_1^{-1}:\Phi X\to\Phi R$. Furthermore, $\Phi q\cdot q^{-1}=1_{\Phi X}=\Phi\pi_1\cdot\pi_1^{-1}$. It is left to show that
\begin{equation*}
 q^{-1}\cdot\Phi q=\Phi \pi_2\cdot\pi_1^{-1},
\end{equation*}
which can be shown with the same calculation as in \cite{Hof_Cocompl}, based on the following proposition.
\begin{proposition}\label{proper}
Consider the following diagram in $\Cat{\Tth}$
\begin{equation*}
\xymatrix{R\ar@<0.5ex>[r]^{\pi_1}\ar@<-0.5ex>[r]_{\pi_2}
&X\ar[r]^{q} &Q}
\end{equation*}
with $\pi_1,\pi_2:R\rightrightarrows X$ in $\CocomplSepPhi{\Tth}$, $(\pi_1, \pi_2)$ an equivalence relation in $\SET$, and $q:X\to Q$ its coequaliser in $\Cat{\Tth}$.
\begin{enumerate}
\item If $\pi_1,\pi_2$ are left adjoints, then $q$ is proper.
\item The diagram
\begin{equation*}
\xymatrix{\Phi R\ar@<0.5ex>[r]^{\Phi\pi_1}\ar@<-0.5ex>[r]_{\Phi\pi_2}
&\Phi X\ar@/_1.5pc/[l]_{\pi_1^{-1}}\ar[r]^{\Phi q} &\Phi Q\ar@/_1.5pc/[l]_{q^{-1}}}
\end{equation*}
is a split fork in $\Cat{\Tth}$.
\end{enumerate}
\end{proposition}
\begin{proof}
(1) As in \cite[Lemma 2.19 and Corollary 2.20]{Hof_Cocompl}.

(2) Analogous to the proof presented in \cite[Section 2.6]{Hof_Cocompl}.
\end{proof}
Finally, we conclude that:
\begin{theorem}\label{monadic}
Under (Ax~1)-(Ax~4), the forgetful functor $G:\CocomplSepPhi{\Tth}\to\SET$ is monadic.
\end{theorem}
\begin{proof}
In order to show that $\CocomplSepPhi{\Tth}$ has and $G$ preserves
coequaliser of $G$-equivalence relations, consider again the first diagram of Proposition \ref{proper}. 
We have seen that
\begin{equation*}
\xymatrix{\Phi R\ar@<0.5ex>[r]^{\Phi\pi_1}\ar@<-0.5ex>[r]_{\Phi\pi_2}
&\Phi X\ar@/_1.5pc/[l]_{\pi_1^{-1}}\ar[r]^{\Phi q} &\Phi Q\ar@/_1.5pc/[l]_{q^{-1}}}
\end{equation*}
is a split fork and hence a coequaliser diagram in $\Cat{\Tth}$. Since
$\pi_1$ and $\pi_2$ are $\Phi$-cocontinuous, there is a $\Tth$-functor
$\Sup_Q^\Phi:\Phi Q\to Q$ which, since $q:X\to Q$ is the coequaliser of
$\pi_1,\pi_2:R\rightrightarrows X$ in $\Cat{\Tth}$, satisfies
$\Sup_Q^\Phi\cdot\yoneda_Q^\Phi=1_Q$. The situation is depicted in the
following diagram.
\[
\xymatrix{R\ar@<0.5ex>[r]^{\pi_1}\ar@<-0.5ex>[r]_{\pi_2}\ar[d]_{\yoneda_R^\Phi}
& X\ar[r]^q\ar[d]^{\yoneda_X^\Phi} &
Q\ar[d]^{\yoneda_Q^\Phi}\ar@/^3pc/[dd]^{1_Q}\\
\Phi R\ar@<0.5ex>[r]^{\Phi \pi_1}\ar@<-0.5ex>[r]_{\Phi
\pi_2}\ar[d]_{\Sup_R^\Phi} &
\Phi X\ar[r]^{\Phi q}\ar[d]^{\Sup_X^\Phi} &\Phi
Q\ar@{..>}[d]^{\Sup_Q^\Phi}\\
R\ar@<0.5ex>[r]^{\pi_1}\ar@<-0.5ex>[r]_{\pi_2} & X\ar[r]^q & Q}
\]
We conclude that $Q$ is separated and $\Phi$-cocomplete, and $q:X\to Q$
is $\Phi$-cocontinuous. Finally, to see that $q:X\to Q$ is the
coequaliser of $\pi_1,\pi_2:R\rightrightarrows X$ in
$\CocomplSepPhi{\Tth}$, let $h:X\to Y$ be in $\CocomplSepPhi{\Tth}$ with
$h\cdot\pi_1=h\cdot\pi_2$. Then, since $\Phi q$ is the coequaliser of
$\Phi\pi_1,\Phi\pi_2:\Phi R\rightrightarrows \Phi X$ in
$\CocomplSepPhi{\Tth}$, there is a $\Phi$-cocontinuous $\Tth$-functor
$f:\Phi Q\to Y$ such that $f\cdot\Phi q=h\cdot \Sup_X^\Phi$. Then
\[
f\cdot\yoneda_Q^\Phi\cdot q=f\cdot\Phi q\cdot\yoneda_X^\Phi=h\cdot
\Sup_X^\Phi\cdot\yoneda_X^\Phi=h
\]
and
\begin{multline*}
\Sup_Y^\Phi\cdot\Phi f\cdot\Phi \yoneda_Q^\Phi\cdot\Phi q
= f\cdot \Sup_{\Phi Q}^\Phi\cdot\Phi \yoneda_Q^\Phi\cdot\Phi q
= f\cdot\Phi q
= h\cdot \Sup_X^\Phi\\
= f\cdot \yoneda_Q^\Phi\cdot q\cdot \Sup_X^\Phi
= f\cdot \yoneda_Q^\Phi\cdot \Sup_Q^\Phi\cdot \Phi q,
\end{multline*}
hence $\Sup_Y\cdot\Phi (f\cdot\yoneda_Q^\Phi)=f\cdot\yoneda_Q^\Phi\cdot
\Sup_Q^\Phi$, that is, $f\cdot\yoneda_Q^\Phi$ is $\Phi$-cocontinuous.
\end{proof}

\section{The examples}

\subsection{All distributors}
The class $\Phi=\calP$ of all distributors satisfies obviously all four axioms. In fact, this is the situation studied in \cite{Hof_Cocompl}.

\subsection{Representable distributors}
The smallest possible choice is $\Phi=\calR$ being the class of all representable $\Tth$-distributors $\calR=\{f^*\mid f\text{ is a $\Tth$-functor}\}$. Clearly, $\calR$ satisfies (Ax~1), (Ax~2) and (Ax~3) but not (Ax~4). We have $\calR(X)=\{x^*\mid x\in X\}$, each $\Tth$-category is $\calR$-cocomplete and each $\Tth$-functor is $\calR$-cocontinuous, and therefore $\CocomplSepVar{\calR}{\Tth}=\CatSep{\Tth}$. This case is certainly not very interesting; however, our results tell us that the inclusion functor $\CatSep{\Tth}\hrw\Cat{\Tth}$ is monadic. In particular, {\em the category $\TOP_0$ of topological T$_0$-spaces and continuous maps is a monadic subcategory of $\TOP$}.

\subsection{Almost representable distributors}
We can modify slightly the example above and consider $\Phi=\calR_0$ the class of all almost representable $\Tth$-distributors, where a $\Tth$-distributor $\varphi:X\kmodto Y$ is called {\em almost representable} whenever, for each $y\in Y$, either $y^*\kleisli\varphi=\bot$ or $y^*\kleisli\varphi=x^*$ for some $x\in X$. As above, $\calR_0$ satisfies (Ax~1), (Ax~2) and (Ax~3) but not (Ax~4).

By definition, for a $\Tth$-category $X$ we have
\[
 \calR_0(X)=\{\psi\in\Psh X\mid \psi\in\calR_0\}=\{x^*\mid x\in X\}\cup\{\bot\},
\]
with the structure inherited from $\Psh X$. Furthermore, a $\Tth$-functor $f:(X,a)\to(Y,b)$ is $\calR_0$-dense whenever, for each $y\in Y$,
\[
\exists\frx\in TX\,.\,b(Tf(\frx),y)>\bot\;\Rw\; \exists x\in X\,\forall\frx\in TX\,.\, b(Tf(\frx),y)=a(\frx,x).
\]
Hence, with
\[
Y_0=\{y\in Y\mid \exists\frx\in TX\,.\,b(Tf(\frx),y)>\bot\}
\]
we can factorise an $\calR_0$-dense $\Tth$-functor $f:X\to Y$ as
\[
X\xrightarrow{\hspace{1ex}f\hspace{1ex}}Y_0\hrw Y,
\]
where $Y_0\hrw Y$ is fully faithful and $X\stackrel{f}{\to} Y_0$ is left adjoint. If we consider $f:X\to Y$ in $\TOP$, then $Y_0=\overline{f(X)}$ is the closure of the image of $f$, so that each $\calR_0$-dense continuous map factors as a left adjoint continuous map followed by a closed embedding. Consequently, {\em for a topological space $X$, the following assertions are equivalent:
\begin{eqcond}
\item $X$ is injective with respect to $\calR_0$-dense fully faithful continuous maps.
\item $X$ is injective with respect to closed embeddings.
\end{eqcond}}
Note that in this example we are working with the dual order, compared with \cite[Section 11]{EF_SemDom}.

\subsection{Right adjoint distributors}
Now we consider $\Phi=\calL$ the class of all right adjoint $\Tth$-distributors. This class contains all distributors of the form $f^*$, for a $\Tth$-functor $f$, and it is closed under composition. Since adjointness of a $\Tth$-distributor $\varphi:X\kmodto Y$ can be tested pointwise in $Y$, the axioms (Ax~1), (Ax~2) and (Ax~3) are satisfied. By definition, $\calL(X)=\{\psi\in\Psh X\mid \psi\text{ is right adjoint}\}$, and a $\Tth$-category is $\calL$-cocomplete if each pair $\varphi\dashv\psi$, $\varphi:Y\kmodto X$, $\psi:X\kmodto Y$, of adjoint $\Tth$-distributors is of the form $f_*\dashv f^*$, for a $\Tth$-functor $f:Y\to X$. For $\V$-categories, this is precisely the well-known notion of Cauchy-completeness as introduced by Lawvere in \cite{Law_MetLogClo} as a generalisation of the classical notion for metric spaces. However, Lawvere never proposed the name ``Cauchy-complete'', and, while working on this notion in the context of $\Tth$-categories in \cite{CH_Compl} and \cite{HT_LCls}, we used instead Lawvere-complete and L-complete, respectively. Furthermore, one easily verifies that each $\Tth$-functor is $\calL$-cocontinuous, i.e.\ (right adjoint)-weighted colimits are absolute, so that $\CocomplSepVar{\calL}{\Tth}=\CatCompl{\Tth}$ is the full subcategory of $\Cat{\Tth}$ consisting of all separated and Lawvere complete $\Tth$-categories.

On the other hand, for a surjective $\Tth$-functor $f$, $f_*$ does not need to be right adjoint, so that (Ax~4) is in general not satisfied. This is not a surprise, since natural instances of this example fail Theorem \ref{monadic}. Indeed, in the category of ordered sets and monotone maps, any ordered set is Lawvere-complete, hence the category of Lawvere-complete and separated ordered sets coincides with the category of anti-symmetric ordered sets. The canonical forgetful functor from this category to $\SET$ is surely not monadic. Also, the canonical forgetful functor from the category of Lawvere-complete and separated topological spaces (= sober spaces) and continuous maps to $\SET$ is also not monadic.

\subsection{Inhabited distributors}\label{ExInh}
Another class of distributors considered in \cite{Hof_Cocompl} is $\Phi=\calI$ the class of all inhabited $\Tth$-distributors. Here a $\Tth$-distributor $\varphi:X\kmodto Y$ is called \emph{inhabited} if
\[
\forall y\in Y\,.\, k\le\bigvee_{\frx\in TX}\varphi(\frx,y).
\]
(Ax~3) is satisfied by definition, and in \cite{Hof_Cocompl} we showed already the validity of (Ax~1) and (Ax~2). Furthermore, one easily verifies that (Ax~4) is satisfied. Hence, as already observed in \cite{Hof_Cocompl}, all results stated in Section \ref{results} are available for this class of distributors. Let us recall that, specialised to $\TOP$, inhabited-dense continuous maps are precisely the topologically dense continuous maps, and the injective spaces with respect to topologically dense embeddings are known as \emph{Scott domains} \cite{Book_ContLat}.

\subsection{``closed'' distributors}
A further interesting class of distributors is given by
\[
\Phi=\{\varphi:X\kmodto Y\mid \forall y\in Y,\,\frx\in TX\,.\, \varphi(\frx,y)\le\bigvee_{x\in X}a(\frx,x)\otimes\varphi(e_X(x),y)\},
\]
that is, $\varphi\in\Phi$ if and only if $\varphi\le \varphi\cdot e_X\cdot a$. Clearly, (Ax~3) is satisfied. For each $\Tth$-functor $g:(Y,b)\to(X,a)$ we have
\[
g^*\cdot e_X\cdot a=g^\circ\cdot a\cdot e_X\cdot a\ge g^\circ\cdot a=g^*,
\]
hence $g^*\in\Phi$. Furthermore, given $\Tth$-distributors $\varphi:X\kmodto Y$ and $\psi:Y\kmodto Z$ in $\Phi$, then
\begin{multline*}
\psi\kleisli\varphi=\psi\cdot\Txi\varphi\cdot m_X^\circ
\le \psi\cdot e_Y\cdot b\cdot\Txi\varphi\cdot m_X^\circ
=\psi\cdot e_Y\cdot\varphi
\le\psi\cdot e_Y\cdot\varphi\cdot e_X\cdot a\\
\le\psi\cdot\Txi\varphi\cdot e_{TX}\cdot e_X\cdot a
\le\psi\cdot\Txi\varphi\cdot m_X^\circ\cdot e_X\cdot a
=(\psi\kleisli\varphi)\cdot e_X\cdot a
\end{multline*}
and therefore also $\psi\kleisli\varphi\in\Phi$. We have seen that this class of distributors satisfies (Ax~1), (Ax~2) and (Ax~3). On the other hand, (Ax~4) is not satisfied.

By definition, a $\Tth$-functor $f:(X,a)\to(Y,b)$ is $\Phi$-dense whenever, for all $\frx\in TX$ and $y\in Y$,
\[
b(Tf(\frx),y)\le\bigvee_{x\in X}a(\frx,x)\otimes b(e_Y(f(x)),y).
\]
Hence, each proper $\Tth$-functor (see \cite{CH_EffDesc}) is $\Phi$-dense. In fact, $\Phi$-dense $\Tth$-functors can be seen as ``proper over $\Cat{\V}$'', and the condition above states exactly properness of $f$ if the underlying $\V$-category $\ForgetToV Y$ of $Y=(Y,b)$ is discrete. Furthermore, each surjective $\Phi$-dense $\Tth$-functor is final with respect to the forgetful functor $\ForgetToV:\Cat{\Tth}\to\Cat{\V}$. To see this, let $f:(X,a)\to(Y,b)$ be a surjective $\Phi$-dense $\Tth$-functor, $Z=(Z,c)$ a $\Tth$-category and $g:\ForgetToV Y\to\ForgetToV Z$ a $\V$-functor such that $gf$ is a $\Tth$-functor. We have to show that $g$ is a $\Tth$-functor. Let $\fry\in TY$ and $y\in Y$. Since $Tf$ is surjective, there is some $\frx\in TX$ with $Tf(\frx)=\fry$. We conclude
\begin{align*}
b(\fry, y) &= b(Tf(\frx),y)\\
&\le\bigvee_{x\in X}a(\frx,x)\otimes b(e_Y(f(x)),y)\\
&\le\bigvee_{x\in X}c(T(gf)(\frx),gf(x)\otimes c(e_Z(gf(x)),g(y))\\
&\le c(Tg(\fry),g(y)).
\end{align*}

\subsection{Further examples}
A wide class of examples of injective topological spaces is described in \cite{EF_SemDom}, where the authors consider injectivity with respect to a class of embeddings $f:X\to Y$ such that the induced frame morphism $f_*:\Omega X\to\Omega Y$ preserves certain suprema. A similar construction can be done in our setting; to do so \emph{we assume from now on $T1=1$}. For a $\Tth$-category $X$, the {\em $\V$-category of covariant presheafs} $\V^X$ is defined as
\[
\V^X=\{\alpha:1\kmodto X\mid \text{$\alpha$ is a $\Tth$-distributor}\}
=\{\alpha:X\to\V\mid \text{$\alpha$ is a $\Tth$-functor}\},
\]
and the $\V$-categorical structure $[\alpha,\beta]\in\V$ is given as the lifting
\[
\xymatrix{X & 1,\ar@{-_{>}}|-{\object@{o}}[l]_\beta \ar@{.^{>}}|-{\object@{o}}[dl]^{\alpha\homkleislileft\beta=:[\alpha,\beta]} \\
1\ar@{-^{>}}|-{\object@{o}}[u]^\alpha}
\]
for all $\alpha,\beta\in\V^X$. Since $e_1:1\to T1$ is an isomorphism, this lifting of $\Tth$-distributors does exist and can be calculated as the corresponding lifting of $\V$-distributors
\[
\xymatrix{X & 1.\ar|-{\object@{o}}[l]_{\beta} \ar@{.{>}}|-{\object@{o}}[dl] \\
1\ar|-{\object@{o}}[u]^\alpha}
\]
Each $\Tth$-distributor $\varphi:X\kmodto Y$ induces a $\V$-functor
\[
 \varphi\kleisli (-):\V^X\to\V^Y,\; \alpha\mapsto\varphi\kleisli\alpha,
\]
which is right adjoint if $\varphi$ is a right adjoint $\Tth$-distributor. Given now a class $\Psi$ of $\V$-distributors, we may consider the class $\Phi$ of all those $\Tth$-distributors $\varphi$ for which $\varphi\kleisli (-)$ preserves $\Psi$-weighted limits. This class of $\Tth$-distributors is certainly closed under composition, and contains all right adjoint $\Tth$-distributors, hence it includes all representable ones. Finally, if $\Psi$-weighted limits are calculated pointwise in $\V^X$, then also (Ax~3) is fulfilled. As particular examples we have the class $\Phi$ of all $\Tth$-distributors $\varphi:X\kmodto Y$ for which $\varphi\kleisli (-)$ preserves
\begin{enumerate}
\item the top element of $\V^X$, that is, for which $\varphi\kleisli\top=\top$. In pointwise notation, this reads as
\[
\forall y\in Y\,.\,\top=\bigvee_{\frx\in TX}\varphi(\frx,y)\otimes\top.
\]
If $k=\top$, then this class of $\Tth$-distributors coincides with the class of inhabited $\Tth$-distributors considered in \ref{ExInh}.
\item cotensors, that is, for each $u\in\V$ and each $\alpha\in\V^X$, $\varphi\kleisli\hom(u,\alpha)=\hom(u,\varphi\kleisli\alpha)$.
\item finite infima (cf. \cite[Section 6]{EF_SemDom}).
\item arbitrary infima (cf. \cite[Section 7]{EF_SemDom}).
\item codirected infima (cf. \cite[Section 8]{EF_SemDom}).
\end{enumerate}

\end{document}